\documentclass[10pt]{amsart}

\usepackage{amsthm,amsmath,amssymb}
\usepackage{datetime}
\usepackage{multirow}
\usepackage{graphicx}
\usepackage{booktabs}
\usepackage{calc}
\usepackage{caption}
\usepackage[numbers]{natbib}

\newcommand\Tstrut{\makebox[0pt][c]{\rule{0pt}{2.6ex}}}         % = `top' strut
\newcommand\Bstrut{\makebox[0pt][c]{\rule[-1.2ex]{0pt}{0pt}}}   % = `bottom' strut

\makeatletter
\patchcmd{\@settitle}{\uppercasenonmath\@title}{}{}{}
\patchcmd{\@setauthors}{\MakeUppercase}{}{}{}
\patchcmd{\section}{\scshape}{}{}{}
\makeatother

\title{The Integer Cuboid Table}
\author{Randall L. Rathbun}
\email{randallrathbun@gmail.com}
\subjclass[2010]{11D09,11Y50,11Y70,14G04}
\keywords{integer cuboid, Euler cuboid, perfect integer cuboid}

\begin{document}

\begin{abstract}
Integer cuboids are rectangular Diophantine parallelepipeds It has been discovered
that these cuboids come in 3 varieties: {\it Euler} or {\it body} type,
{\it edge} type, and {\it face} type. In all three cases, one edge or diagonal
is irrational, all others are rational. We discuss an exhaustive computer search procedure which uses the Pythagorean group $\mathnormal{Py(n)}$ to locate all possible cuboids with a given edge $n$. Over the range of 44 to 200,000,000,027 for the smallest edge, 167,043 cuboids were discovered. They are listed in the Integer Cuboid Table.
\end{abstract}

\maketitle

\setcounter{section}{0}
\section*{Diophantine Parallelepipeds}

Parallelepipeds have been examined, in Diophantine analysis of Number Theory.

We define a rational parallelepiped in $n$-dimensions as a polytope spanned by $n$ vectors
$\boldsymbol{\vec{v}}_1, \dots,\:\boldsymbol{\vec{v}}_n$ in a vector space over the rationals,
$\mathbb{Q}$, or integers, $\mathbb{Z}$
\begin{center}
   span($\boldsymbol{\vec{v}}_1,\dots,\:\boldsymbol{\vec{v}}_n) \: = \: t_1\boldsymbol{\vec{v}}_1 + \: \dots + t_n\boldsymbol{\vec{v}}_n \text{ for } t_i \in \mathbb{Q}$ for $i=1 \to n$
\end{center}

Here we are interested in three dimensions, $n=3$, so a rational parallelepiped is a prism determined by the 3 basis vectors
$\boldsymbol{\vec{a}}$, $\boldsymbol{\vec{b}}$, $\boldsymbol{\vec{c}}$. The prism has 8 vertices, 3 pairs of parallel faces
which are all parallelograms, and 12 edges $\in\mathbb{Q}$.
\begin{equation*}
   \text{span}(\boldsymbol{\vec{a}},\boldsymbol{\vec{b}},\boldsymbol{\vec{c}}) \: = \: t_1\boldsymbol{\vec{a}} + t_2\boldsymbol{\vec{b}} + t_3\boldsymbol{\vec{c}} \;\; \text{ for } t_{i=1 \to 3} \in \mathbb{Q}
\label{eq:span}
\end{equation*}
\begin{figure}[!ht]
\centering
\includegraphics[scale=1.0]{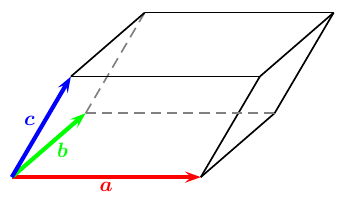}
\captionsetup{justification=centering}
\caption[caption]{The rational parallelepiped with 3 basis vectors $\vec{a}$, $\vec{b}$, and $\vec{c}$}
\label{fig:ratpiped}
\end{figure}
In Figure \ref{fig:ratpiped}, the lengths of the basis vectors, $\vec{a}$, $\vec{b}$, and $\vec{c}$, are rational
$\in \mathbb{Q}$, and they determine the parallelepiped, thus all 12 edges are rational.

For pipeds in general, there are 3 distinct edge lengths, 6 distinct face diagonals, and 4 distinct body diagonals.

What is not so well known, is that there are 5 classes of Diophantine parallelepipeds, based upon the number of right angles
which occur at the intersection of the 3 basis vectors. The classes are triclinic-obtuse, triclinic-acute,
bi-clinic 1-orthogonal, monoclinic 2-orthogonal, and the well known rectangular cuboid which is tri-orthogonal.

\section*{Introduction to integer cuboids}

The {\it integer cuboid} is actually a Diophantine piped that has three right-angles at the intersection of the basis
vectors, $\boldsymbol{\vec{a}}$, $\boldsymbol{\vec{b}}$, $\boldsymbol{\vec{c}}$. It is tri-orthogonal. The piped figure
is rectangular, like a brick, or {\it cuboid}, hence the name.

Due to the three-fold orthogonality, the face diagonals in the piped change from 6 distinct lengths in the
parallelogram faces to just 3 lengths in the rectangle faces. The 4 body diagonals become just 1 body diagonal.
There are still 3 distinct edges to consider.

\begin{figure}[ht]
\centering
\includegraphics[scale=1.15]{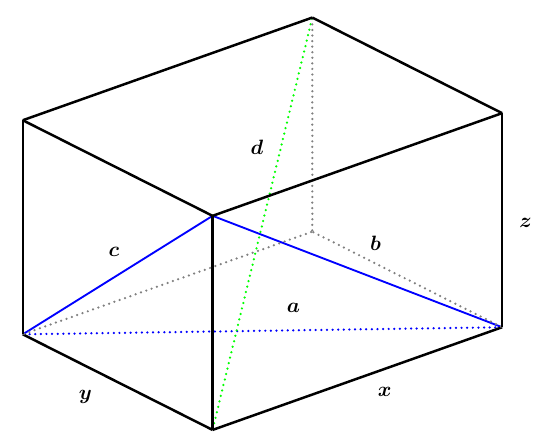}
\caption{The Integer Cuboid.}
\label{fig:cuboid}
\end{figure} \hfill

Thus there are 7 lengths involved in the Diophantine rectangular integer cuboid.
They are the 3 edges: $x$, $y$, and $z$, the 3 face diagonals: $a$, $b$, and $c$ and the body diagonal: $d$.
See figure \ref{fig:cuboid}.

These lengths satisfy the following equations:
\begin{align}
x^2 + y^2 & = a^2 \label{eq:xyz} \\
x^2 + z^2 & = b^2 \label{eq:xzb} \\
y^2 + z^2 & = c^2 \label{eq:yzc} \\
x^2 + y^2 + z^2 & = d^2 \label{eq:diag}
\end{align}
\indent The last equation (\ref{eq:diag}) also means that
\begin{equation}
x^2 + c^2 = y^2 + b^2 = z^2 + a^2 = d^2
\end{equation}
\indent It has been discovered that there are 3 types of cuboids, called {\it body}, {\it edge}, and {\it face} cuboids.
Each type has one of the 7 lengths irrational, the other 6 lengths are rational.

In the case of the body cuboid, the body diagonal $d$ is irrational. For the edge cuboid,
one of the edges $x,y,z$ is irrational. The face cuboid has just one of the face diagonals $a,b,c$ irrational.

The body cuboid is commonly referred to as the {\it Euler cuboid} in honor of Leonard Euler, who discussed this type
of cuboid\cite{euler1}. He was also aware of face cuboids, and provided the 104, 153, 672, 697 example\cite{euler2}.

Only recently have edge cuboids become known\cite{ama,brom,leech}.

A {\it perfect integer cuboid} would have all 7 lengths rational. It is still unknown if it exists or not.

\section*{Setting up the Pythagorean group $\mathnormal{Py(n)}$}

How can we discover integer cuboids? Manual searches have been done, Maurice Kraitchik's Rational Cuboid
Table\cite{krait1,krait2,krait3} is a classic example which took thousands of hours to produce. His table
has been corrected and expanded with the help of a computer\cite{rathbun}. So we will enlist the aid
of a computer to help in expediting that search process.

First of all, we take advantage of an equation to help find pairs of numbers such that when their squares is added,
their sum is a square. The equation is very similar to the Pythagorean equation:

For any $n,\:d \in \mathbb{Q}$, let:
\begin{equation}
m = \frac{n^2-d^2}{2d}
\label{eq:kndiv}
\end{equation}
then $m\in \mathbb{Q}$ and
\begin{align}
m^2 + n^2 & = \left(\frac{n^2-d^2}{2d}\right)^2 + n^2 \\
m^2 + n^2 & = \frac{\left(n^2-d^2\right)^2}{4d^2} + n^2 \\
m^2 + n^2 & = \frac{n^4-2d^2n^2+d^4 + 4d^2n^2}{4d^2} \\
m^2 + n^2 & = \frac{\left(n^2+d^2\right)^2}{4d^2} \\
m^2 + n^2 & = \left(\frac{n^2+d^2}{2d}\right)^2 = \square
\end{align}
If $m \in \mathbb{Q}$ is determined by eq(\ref{eq:kndiv}), using $n,d$, then $m^2+n^2=\square$.

Secondly, we create a set of divisors $d(n)=d_{i=1\to k}$ for any integer $n \in \mathbb{Z}$ according to eqs (\ref{eq:divs1},\ref{eq:divs2}).

\noindent Let
\begin{align}
d(n) & = d_{i=1\to k} = \{ d \: | \: n^2 \} \; \text{:} \; n\bmod 2 = 1 \label{eq:divs1} \\
d(n) & = d_{i=1\to k} = 2\cdot \{ d \:| \left( n/2 \right)^2 \: \} \; \text{:} \; n\bmod 2 = 0 \label{eq:divs2}
\end{align}
where $d(n) = d_{i=1\to k}$ is the set of proper divisors of either $n^2$ when $n$ is odd or twice the set of proper divisors of
$\left(\frac{n}{2}\right)^2$ when $n$ is even, and $k$ is the count of the proper divisors unique for each $n \in \mathbb{Z}^+$.

Example:
\begin{align*}
n & = 44 \;\: \text{ then } d(44) \;\: = \{ 2, 4, 8, 22 \} \quad d_{i=1\to 4} \\
n & = 117 \text { then } d(117) = \{ 1, 3, 9, 13, 27, 39, 81 \} \quad d_{i=1\to 7}
\end{align*}

Having found the set $d(n) = d_{i=1\to k}$ of proper divisors of any integer $n$, we create a set called
the {\it Pythagorean group} $\mathnormal{Py(n)}$. Using eq(\ref{eq:divs1}), or eq(\ref{eq:divs2}), and letting
$d = d_{i=1\to k}$ for all proper divisors of $n^2$ as defined by eq(\ref{eq:kndiv}), we create 2 sets of numbers:
\begin{align}
a_{i=1\to k} & = \frac{n^2-d_i^2}{2d_i} \label{eq:a} \\
A_{i=1\to k} & = \frac{n^2+d_i^2}{2d_i} \label{eq:A}
\end{align}
then the Pythagorean group is
\begin{equation}
\mathnormal{Py(n)} = \{ a_i, A_i \} \quad i=1\to k
\end{equation}
for $k$ divisors of $n^2$ or $(n/2)^2$, where $a_i,A_i$ are found by eqs (\ref{eq:divs1} or \ref{eq:divs2}, \ref{eq:a}, \ref{eq:A}).

\noindent Using the example above, we have:

$\mathnormal{Py}(44)  = \{ a_i,A_i \}_{i=1\to 4} = \{ [ 483, 240, 117, 33 ], [ 485, 244, 125, 55 ] \}$

$\mathnormal{Py}(117) = \{ a_i,A_i \}_{i=1\to 7} = \{ [ 6844, 2280, 756, 520, 240, 156, 44 ],$ \hfill \\
\phantom{000000000000000000000000000}$\:[ 6845, 2283, 765, 533, 267, 195, 125 ] \}$ \hfill \\
We also note that
\begin{equation}
a_i + d_i = A_i \; \text{for} \; i=1\to k
\end{equation}
so $A_i$ does not have to be explicitly derived as in eq(\ref{eq:A}), but found from $a_i + d_i$.

\section*{Searching for Cuboids}

We start with an edge $N$ under computer consideration, and find its Pythagorean group
$\mathnormal{Py(N)}$.

\noindent Let
\begin{equation}
\mathnormal{Py(N)} = \{ a_i, A_i \} \quad i=1\to k
\end{equation}
be the Pythagorean group for the edge. By examining certain square sums among selected 
$A_{i,j}\;a_{i,j} : \:\scriptstyle{(i<j\le k)}$ pairs of the group which satisfy the conditions as given in
the search condition table \ref{table:search} below, we can discover cuboids which may exist.

\begin{table}[ht]
\centering
Condition table for possible cuboids in the Pythagorean group \\[0.5em]
\begin{tabular}{c|c|c}
condition & cuboid solution & type \\
\toprule
  $A_i^2 + a_j^2 = s^2$  &  $a_i,\;a_j,\;N,\;s$                &  face\Tstrut\Bstrut \\
\midrule
  $A_i^2 - A_j^2 = s^2$  &  $s,\;N,\;a_j,\;A_i$                &  face\Tstrut\Bstrut \\
\midrule
  $A_j^2 - A_i^2 = s^2$  &  $s,\;N,\;a_i,\;A_j$                &  face\Tstrut\Bstrut \\
\midrule
  $A_i^2 - a_j^2 = s^2$  &  $N,\;a_j,\;\sqrt{s^2-N^2},\;A_i$   &  edge\Tstrut\Bstrut \\
\midrule
  $a_j^2 - A_i^2 = s^2$  &  $N,\;a_j,\;\sqrt{-s^2-N^2},\;A_i$  &  edge\Tstrut\Bstrut \\
\midrule
  $A_j^2 - a_i^2 = s^2$  &  $N,\;a_i,\;\sqrt{s^2-N^2},\;A_j$   &  edge\Tstrut\Bstrut \\
\midrule
  $a_i^2 - A_j^2 = s^2$  &  $N,\;a_i,\;\sqrt{-s^2-N^2},\;A_j$  &  edge\Tstrut\Bstrut \\
\midrule
  $a_i^2 + a_j^2 = s^2$  &  $a_i,\;a_j,\;N,\;\sqrt{s^2+N^2}$   &  body\Tstrut\Bstrut \\
\bottomrule
\end{tabular}\vspace{0.1in}
\caption{Cuboid Search Table}
\label{table:search}
\end{table}

The cuboid solution is given as three edges $x$, $y$, and $z$, and the body diagonal $d$, in
the cuboid solution column. Please note that $s\in\mathbb{Z}^{+}$ is integer. If the value inside
the radical is negative, the we have an edge cuboid with one edge which is a complex number,
otherwise the edge cuboid is real.

For each edge $N$, the computer has to create the Pythagorean group $\mathnormal{Py(N)}$, then
examine all possible pairs $A_{i,j}$, and $a_{i,j}$ for the conditions listed above. This is a
type $\mathit{O}^2$ search because all possible $i,j$ pairs needs to be checked. Using modulo
conditions upon the $A_{i,j}$ or $a_{i,j}$ pairs considerably shortens down the testing.

It has to be noted that once a cuboid is found, it is reduced to primitive terms, taking care to
appropriately reduce the value inside the radical, if the cuboid is an edge or Euler(body) type.

\noindent Example:
\begin{equation*}
\mathnormal{Py}(44) = \{ a_i,A_i \}_{i=1\to 4} = \{ [ 483, 240, 117, 33 ], [ 485, 244, 125, 55 ] \} \\
\end{equation*}
The computer discovers that $240^2+117^2=267^2$ for $a_2$ and $a_3$ and we have the last condition satisfied,
so we discover the Euler(body) cuboid 44, 117, 240, $\sqrt{73225}$.

\noindent Another example:
\begin{multline*}
\mathnormal{Py}(104) = \{ a_i,A_i \}_{i=1\to 7} = \{ [ 2703, 1350, 672, 330, 195, 153, 78 ],\\
 [ 2705, 1354, 680, 346, 221, 185, 130 ] \} \;\; \quad \quad \quad \quad
\end{multline*}
The computer discovers that $680^2+153^2=697^2$ for $A_3$ and $a_6$ and we have the first condition satisfied,
and the face cuboid $153,\; 672,\; 104,\; 697$ is found. One face diagonal is irrational, $\sqrt{474993}$.

Similarly, the computer exhaustively selects each edge $N$, finds its Pythagorean group $\mathnormal{Py(N)}$, examines the
$A_{i,j}$ and $a_{i,j}$ pairs and collects each primitive cuboid which satisfies any 1 of the 8 conditions given in the
search condition table \ref{table:search} above. The use of parallel processors greatly speeds up this task.

\section*{The Integer Cuboid Table}

The {\it Integer Cuboid Table} is a list of all primitive cuboids found by exhaustive computer
search, using the Pythagorean group $\mathnormal{Py(N)}$ and the search condition table \ref{table:search},
for all integers $N$ from 44 to 200,000,000,000. The actual table accompanies this introductory paper as the gzipped
archive {\it integer\_cuboid\_table.gz}.

\subsection*{The Integer Cuboid Table Format}

The table is stored as tab-delimited ASCII characters, as shown below in Table 2.
\begin{table}[!ht]
\centering
\caption{Integer cuboid table - raw ascii data format}
\begin{tabular}{rrl}
\# & ss & cuboid \\
1 & 44 & B,44,117,240,(73225) \\
2 & 60 & e,60,63,(-3344),65 \\
3 & 85 & B,85,132,720,(543049) \\
4 & 104 & F,153,672,104,697 \\
5 & 108 & e,108,725,(-426400),333 \\
6 & 117 & F,520,756,117,925 \\
7 & 124 & E,124,957,(13852800),3845 \\
%8 & 140 & B,140,480,693,(730249) \\
%9 & 160 & B,160,231,792,(706225) \\
%10 & 187 & B,187,1020,1584,(3584425) \\
%11 & 195 & B,195,748,6336,(40742425) \\
%12 & 215 & E,215,912,(533332800),23113 \\
%13 & 240 & B,240,252,275,(196729) \\
%14 & 240 & e,240,364,(-56871),365 \\
%15 & 252 & F,2261,2640,252,3485 \\
%16 & 264 & F,264,952,495,1105 \\
%17 & 264 & F,264,975,448,1105 \\
%18 & 324 & E,324,4368,(23914825),6565 \\
%19 & 333 & F,644,2040,333,2165 \\
\end{tabular}
\end{table}
The number inside a parenthesis pair denotes taking the square root value. Thus (73225) means $\sqrt{73225}$.

\noindent The following detail explains the layout of the Integer Cuboid Table. See Table 3.

\begin{table}[!ht]
\centering
\caption{Explanation of raw data}
\begin{tabular}{r|r|c|rrrr}
\multirow{2}{*}{\#} & \multirow{2}{*}{ss} & \multirow{2}{*}{type} & \multicolumn{4}{c}{cuboid} \\
 & & & x & y & z & d \\
\midrule[0.11em]
1 & 44 & B & 44 & 117 & 240 & \Tstrut $\sqrt{73225}$ \\
2 & 60 & e & 60 & 63 &$\sqrt{-3344}$ & 65 \\
3 & 85 & B & 85 & 132 & 720 & $\sqrt{543049}$ \\
4 & 104 & F & 153 & 672 & 104 & 697 \\
5 & 108 & e & 108 & 725 & $\sqrt{-426400}$ & 333 \\
6 & 117 & F & 520 & 756 & 117 & 925 \\
7 & 124 & E & 124 & 957 & $\sqrt{13852800}$ & 3845 \\
%. & & & & & & \\
%154569 & 154991590800 & F & 154991590800 & 200271206124 & 823272908843 & 861341530205 \\
%154570 & 154994410600 &	B & 154994410600 & 197727052545 & 402857641824 & $\sqrt{225413534201368857384001}$ \\
%154571 & 154995920040 &	F & 154995920040 & 1071739709964 & 194510345323 & 1100219985085 \\
%\bottomrule
\end{tabular}
\
\end{table}

The first column is the index \# for the table. 

The second column 'ss' is the sorted side which has the key value by which the Integer Cuboid Table is sorted.
Some sides are repeated, for cuboids which have the same smallest edge, the sorting was continued by a
numerical sort upon the other edges of these sets of cuboids.

The third column labels the cuboid type. It uses: B for body or Euler cuboid, 
e for an edge cuboid with a complex length, E for a normal edge cuboid, and F for the face cuboid.
The edge cuboid has a irrational length which is a complex number. These solutions were included into
the cuboid table in order to preserve the Pythagorean generator relationships with other cuboids.
The normal E edge cuboid has a real irrational edge.

The next 3 columns are the edge lengths for $x$, $y$, and $z$ from eqs (\ref{eq:xyz},\ref{eq:xzb},\ref{eq:yzc}).
The last column is the length $d$ of the body diagonal, eq (\ref{eq:diag}).

\subsection*{Comments}

The {\it Integer Cuboid Table} is believed complete for this range. Recently another type of cuboid search\cite{DArox}  was implemented
and over the same range, the cuboid counts matched, lending confidence to the completeness of this table.

167,043 cuboids were found: 61,042 were Euler(body) cuboids, 16,612 were edge cuboids with a complex number edge length, 32,286 were edge cuboids,  and 57,103 were face cuboids. The approximate ratios of occurrence was 201 : 161 : 188 for body : edge : face cuboids.

All discovered Euler (body) cuboids have an irrational body diagonal.

A {\it perfect integer cuboid} would have a rational body diagonal and three rational face diagonals. None was found over the range of this search.

The {\it Integer Cuboid Table} for the smallest edge $n$ was exhaustively complete as of Thursday June 25 00:06:37 PDT 2020 over the range $ 1 \le n \le 200,000,000,027$.

\noindent\rule{2.55in}{0.4pt} \vspace*{-0.7em} \par
\settimeformat{hhmmsstime}
{\footnotesize
\noindent Typeset - \currenttime { PDT } \today \par
}


\begin{thebibliography}{10}

	\bibitem{ama} J. Peacock, J. Hancock, N. A. Phillips, {\em Mahatma's problem No. 78},
	Journal Assist. Masters Assoc. London, vol 44 (1949) p118, p225

	\bibitem{brom} Bromhead, T. B., {\em 2918. On square sums of squares},
	Mathematical Notes 2918, Mathematical Gazette, 1960, Vol. 44, Num 349, Oct. pp 219-220

	\bibitem{euler1} Euler, Leonard, {\em Vollst\"{a}ndige Anleitung zur Algebra},
	Kayserliche Akademie der Wissenschaften, St. Petersburg, 1771

	\bibitem{euler2} Euler, Leonard, {\em Vollst\"{a}ndige Anleitung zur Algebra},
	2, Part II, §236, English translation: Euler, Elements of Algebra, Springer-Verlag 1984

	\bibitem{krait1} Kraitchik, Maurice, {\em Th\`{e}orie des Nombres},
	Tome 3, Analyse Diophantine et application aux cuboides rationelles,
	Gauthier-Villars, Paris, 1947

	\bibitem{krait2} Kraitchik, Maurice, {\em On certain rational cuboids},
	Scripta Mathematica, Vol. 11, 1945, pp 317-326

	\bibitem{krait3} Kraitchik, Maurice, {\em Sur les cuboides rationelles},
	Proceedings, International Congress Mathematics, 1954, Vol. 2 Amsterdam, pp 33-34

	\bibitem{rathbun} Rathbun, Randall, {\em The Rational Cuboid Table of Maurice Kraitchik}
	(extended to odd side less than $2^{32}$), arXiv:math/0111229 [math.HO]

	\bibitem{leech} Leech, John, {\em The rational cuboid revisited},
	American Mathematical Monthly, 1977, Vol. 84, pp 518-553,
	Erratum, Amer. Math. Monthly, 1978, Vol. 85, p. 472

	\bibitem{DArox} DArox, Renyxa, {\em Integer cuboid searches}, private communications, July 4, 2017

\end{thebibliography}
\end{document}